\documentclass[runningheads,a4paper]{llncs}

\usepackage{amssymb}
\usepackage{graphicx}

\usepackage{url}
\urldef{\mailsa}\path|jaimecs@mat.uc.pt|    
\newcommand{\keywords}[1]{\par\addvspace\baselineskip
\noindent\keywordname\enspace\ignorespaces#1}

\begin{document}

\mainmatter
\title{What international studies say about \\the importance and limitations of using computers to teach mathematics\\ in secondary schools \thanks{The final publication is available at http://link.springer.com.}}

\titlerunning{Importance and Limitations of Using Computers to Teach Mathematics}
\author{Jaime Carvalho e Silva}
\authorrunning{Conferences on Intelligent Computer Mathematics - CICM 2014}

\institute{Departamento de Matem\'atica, Universidade de Coimbra, \\
Apartado 3008, EC Santa Cruz, 3001-501 Coimbra, Portugal\\
\mailsa\\
\url{http://www.mat.uc.pt/~jaimecs/pessoal/}}

\toctitle{Lecture Notes in Computer Science}

\maketitle

\begin{abstract}
The use of technology in schools has been one of the most debated topics around mathematics education. In some countries there is a huge investment, in others there is a downscaling. Malaysia decided in 2013 to put its 10 million students to use Google laptops and Google apps, while Australia in the same year decided it would not continue funding their own high school laptop program. Who is right from the educational point of view?
The last major curriculum document written in the world to date, the Common Core State Standards-CCSS in the United States, whose mathematics part is coordinated by the well known mathematician William McCallum, sets as one of its standards for mathematical practice:
''Mathematically proficient students consider the available tools when solving a mathematical problem. These tools might include pencil and paper, concrete models, a ruler, a protractor, a calculator, a spreadsheet, a computer algebra system, a statistical package, or dynamic geometry software.''
Strong moves need substantiation from research, including the analysis of the existing situation in different countries. What does research say about the use of computers in schools in present time and the use of different pieces of software from spreadsheets to computer algebra systems?
\keywords{mathematics education}
\end{abstract}

\section{Introduction}

It is extremely frequent to see heated debates on the pros and cons of using computers in the classroom, be it the primary school classroom, the middle school classroom, the secondary school classroom or even the higher education classroom.\\

For example \cite{stross} in the \emph{The New York Times}, on June 10, 2010, Randall Stross, a professor of business at San Jose State University, quotes some studies about the impact of computers and technology on student's test scores to conclude that  

\begin{quote}
Middle school students are champion time-wasters. And the personal computer may be the ultimate time-wasting appliance. Put the two together at home, without hovering supervision, and logic suggests that you wonÕt witness a miraculous educational transformation. 
\end{quote}

One of the studies he quotes was published in the \emph{Quarterly Journal of Economics} and carried out in Romania, comparing families that in 2009 received vouchers to buy computers against families that did not receive that voucher. Results showed that students in the first group showed \emph{significantly lower school grades in math, English and Romanian}. This is a good example of bad use of statistics. You carry out a large scale collection of data, isolating a few variables and ignoring most of them and then imply there is a cause-effect relation between the variables you isolated. Anybody that works closely within education knows you cannot draw simple conclusions from complex data. Another study of a similar kind arrives at similar conclusions: a working paper published by the \emph{National Bureau of Economic Research} correlates data relative to the introduction of broadband services in the american state of North Carolina between 2000 and 2005 to the results of middle school test scores during that period. This time they conclude that the ''negative effects'' were limited to low income households.\\

In a recent analysis \cite{rot} published also in \emph{The New York Times}, Carlo Rotella, director of American studies at Boston College, discusses in detail the position of Joel Klein, chief executive of an IT company and former  chancellor of New York CityÕs public schools from 2002 to 2011. In short, Joel Klein believes that teachers and students need new and interesting tools that help them teach and learn. When asked about evidence for his claims he just says that tablets ''will help teachers do'' what educational research shows is important, that  ''an individual student will learn more if you can tailor the curriculum to match her learning style, pace and interests.'' Of course, he did not prove that tablets will indeed accomplish this goal.\\

Lots of countries are investing hard in introducing some kind of high tech tools in the classroom. John Koetsier reports \cite{koet}, on the digital news site \emph{VentureBeat}, that the Malaysian government is investing massively in introducing computers, internet access software and ''Google Apps'' in all 10 thousand schools in the country totaling 10 million students, teachers, and parents. Why was this choice made? Because ''Google Apps'' are, for educational use, completely free. Is this the best software for educational use? Will students learn better mathematics with this environment? This is not stated.\\

It is clear that we have lots of options on the use of IT technology in schools (including not using it at all), but we need clear ideas before accepting or rejecting hardware, software and communications in the classroom.

\section{International studies}

We live in a society where technology is in a rapid evolution, with new tools arriving at the consumer market every year. It is more than natural that these tools are also offered at the school level. Two main reasons can be stated in favor of this: first of all the school has never been as efficient as the society desires and so new approaches are normally welcome (at least by most people); secondly, if the school is to prepare students for ''real life'' and for some professional activity, then teaching should somehow incorporate the technological tools that students will find someday in their adult life.\\

In education, things are never simple. Several authors, like Luc Trouche \cite{trou}, already pointed out we should make a distinction between an artifact and an instrument. This means that you may have some artifact in the classroom (like a technological tool, hardware or software) but it may have no effect at all unless you are able to integrate it in your activity, and then it becomes an instrument for you. This is not a useless distinction because we have important examples of this difference and the PISA OECD studies give us one of these.\\

\subsection{Digital reading}

PISA is a program conducted by the OECD to study the extent to which 15-year-old students (normally near the end of compulsory education) have acquired the knowledge and skills that are essential for full participation in modern society, focusing mainly in mathematics, reading and science. This program began in the year 2000 and is applied in numerous countries every three years. A lot of data is collected about the students, the teachers, the schools and the student's environment. Also some other studies are conducted in parallel, at least in some of the countries participating in the main PISA study. In 2012 a total of 65 countries and economies participated in the PISA data collecting but only 44 countries and economies participated in a computer-based assessment of problem solving; 32 of them also participated in a computer-based assessment of reading and mathematics.\\

For the first time, the PISA 2009 survey also assessed 15-year-old studentsÕ ability to read, understand and apply digital texts. These texts are very different from printed ones, namely at the level of their organisation. In 19 countries and economies students were given questions via computer to assess this ability. The PISA 2009 results \cite{pisa2011} about digital reading show something striking. They show that even when guidance on navigation is explicit, significant numbers of students still cannot locate crucial pages.\\

Digital reading poses new problems to users: indexing and retrieval techniques are new because of the virtual nature of page contents and formats; also hyperlinks are introduced and new multipage documents are used in a networked structure that may confuse the reader. So, a completely new environment comes up and PISA digital reading assessment offers powerful evidence that todayÕs 15-year-olds, the ''digital natives'', do not automatically know how to operate effectively in the digital environment, contrarily to what we could have thought.\\

\subsection{Computer-based assessment of mathematics}

For the first time in 2012, PISA included an optional computer-based assessment of mathematics. 32 of the 65 countries and economies participated in this. Specially designed PISA questions were presented on a computer, and students responded on the computer, although they could also use pencil and paper as they worked out through the test questions. The PISA 2012 report \cite{pisa2012} justifies this part of the PISA program:

\begin{quote}
(...) computer-based items can be more interactive, authentic and engaging than paper-based items. They can be presented in new formats (e.g. drag-and-drop), include real-world data (such as a large, sortable dataset), and use colour, graphics and movement to aid comprehension. Students may be presented with a moving stimulus or representations of three-dimensional objects that can be rotated, or have more flexible access to relevant information. New item formats can expand response types beyond verbal and written, giving a more rounded picture of mathematical literacy. (...) computers have become essential tools for representing, visualising, exploring, and experimenting with all kinds of mathematical objects, phenomena and processes, not to mention for realising all types of computations Ð at home, at school, and at work.\cite{pisa2012}
\end{quote}

Fourty one specially designed computer-based items were developed for this assessment. These items were designed so that  mathematical reasoning and processes would take precedence over the ability of using the computer as a tool. The report details the approach used:

\begin{quote}
Each computer-based item involves three aspects:\\
- the mathematical demand (as for paper-based items);\\
- the general knowledge and skills related to information and communication technologies (ICT) that are required
(e.g. using keyboard and mouse, and knowing common conventions, such as arrows to move forward). These are
intentionally kept to a minimum;\\
- competencies related to the interaction of mathematics and ICT, such as making a pie chart from data using a simple
''wizard'', or planning and implementing a sorting strategy to locate and collect desired data in a spreadsheet.\cite{pisa2012}
\end{quote}

The conclusion of this part of the study is that ''there is a high degree of consistency in student performance on items delivered on paper and by computer'' but with some important exceptions:

\begin{quote}
In the field of mathematics, one participant (Shanghai-China) saw a large difference, of around 50 score points, in favour of the paper based format. Three other countries and economies showed substantial differences in the same direction - Poland (28-point difference), Chinese Taipei (22-point difference) and Israel (20-point difference). Conversely, there are also countries for which computer delivery of the assessment appears to have been advantageous. The largest difference, of about 30 score points, was seen in Brazil. Colombia also saw a difference of about 20 points in the same direction. The United States, the Slovak Republic and Italy also saw marked, albeit smaller, differences in favour of the computer delivery of the assessment. Across OECD countries, the performance advantage of the computer-based assessment is slightly higher for boys than for girls. (\cite{pisa2012}, p. 491)
\end{quote}

This is a quite recent report and these differences are not yet discussed in terms of the nature of the tasks, of the mode of delivery, or of the student familiarity with computers. In the PISA 2015 program, the computer-based assessment will be the primary mode of delivery for mathematics literacy and all the other domains, but the use of paper-based assessment instruments is an option for countries choosing to do so. In some years we will have then more data for our discussion.

\subsection{Improvements in performance}

Another important question we need to answer, and is raised by a lot of people, is wether students perform better or worse in a computer environment (at school and at home). Some PISA studies also address this question. The PISA 2003 study discusses the relation between the frequency of use of computers at home and student performance in mathematics. And the conclusion is very clear:

\begin{quote}
(...) in every country, students reporting rare or no use of computers at home (on average 18\% of students) score much lower than their counterparts reporting moderate use or frequent use. \cite{pisa2006}
\end{quote}

The PISA 2006 study compares the PISA scores and the use of ICT. In \cite{pisa2010} students are grouped according to frequency of ICT use and then the average performances of each group are compared. Of course this does not tell the whole story because some factors that affect computer use also affect student performance. In order to give a clear picture the PISA 2006 study includes questions about the location and frequency of student computer use. Trying to include a number of relevant variables the PISA 2006 study concludes that:

\begin{quote}
A higher frequency of computer use is associated with higher average science scores in all countries considered. Among OECD countries, the largest effect of using a computer almost every day was found in Iceland, Japan, The Netherlands, Norway, Poland and Spain. Among partner countries, the largest effect of using computer almost every day was found in Bulgaria; Macao, China; and Slovenia. (\cite{pisa2010}, p. 150) (...) in a large majority of countries, the benefits from higher computer use tend to be greater at home than at school. Therefore, despite the better environment and support that schools are expected to provide, computer use tends to have less impact at school than at home. (\cite{pisa2010}, p. 156)
\end{quote}

Having this in mind the PISA study recommends concrete actions regrading ICT use in schools:

\begin{quote}
(...) the analysis has shown that computer use increases student performance but that this increase is not the same for all students. (...) as the benefits from computer use depend on the characteristics of each student, policies to increase ICT use need to be tailored to students. (...) the positive effects of computer use on student performance are greatest when they are supported by a sufficient level of capital. Skills, interests and attitudes affect studentsÕ engagement with ICT, the activities they carry out on the computer and how well. An increase in ICT use that is not supported by an increase in capital would have a lower impact on student performance. (\cite{pisa2010}, p. 156)
\end{quote}

Of course, it is clear from these studies that the simple use of ICT does not guarantee an improvement in performance:

\begin{quote}
(...) the apparently negative association between performance and some kinds of computer usage, shown by PISA 2003 and now PISA 2006, carries a warning not to assume that more is better for studentsÕ performance. (\cite{pisa2010}, p. 158)
\end{quote}

It is clear form these studies that the use of ICT has generally very positive effects on student performance at the mathematics and science level. Only incomplete studies will conclude that the use of ICT has a negative influence in student performance. These OECD big scale studies examine the educational situation in great detail and include very different political and social realities in the big number of countries involved so that their conclusions are very reliable. What we loose in these huge statistical studies is the detail. We need now to know what works and what does not work in each situation.

\section{ICMI studies}

ICMI, the International Commission on Mathematical Instruction, founded in 1908 to foster efforts to improve the quality of mathematics teaching and learning worldwide, has produced two large studies that discuss in detail the impact and use of ICT in mathematics education. These were:

\begin{itemize}
\item ICMI Study 1. The Influence of Computers and Informatics on Mathematics and its Teaching\\
Study Conference held in Strasbourg, France, March 1985.\\
Study Volume published by Cambridge University Press, 1986, eds: R.F. Churchhouse et al. (ICMI Study Series)\\
Second edition published by UNESCO, 1992, eds: Bernard Cornu and Anthony Ralston. (Science and Technology Education No. 44) 
\item ICMI Study 17. Digital Technologies and Mathematics Teaching and Learning: Rethinking the Terrain\\
Study Conference held in Hanoi, Vietnam, December 2006.\\
Study Volume published by Springer, 2010: Mathematics Education and Technology-Rethinking the Terrain. The 17th ICMI Study Series: New ICMI Study Series, Vol. 13. Hoyles, Celia; Lagrange, Jean-Baptiste (Eds.) (New ICMI Study Series 13) 
\end{itemize}

These studies point out some directions for the integration of ICT in mathematics education, but it is also clear that much more research needs to be done:

\begin{quote}
The way digital technologies can support and foster today collaborative work, at the distance or not, between students or between teachers, and also between teachers and researchers, and the consequences that this can have on studentsÕ learning processes, on the evolution of teachersÕ practices is certainly one essential technological evolution that educational research has to systematically explore in the future. (\cite{icmi}, p. 473)
\end{quote}

Numerous examples are described and quoted in this 500-page volume but we need to have in mind what I consider to be the main conclusion:

\begin{quote}
Making technology legitimate and mathematically useful requires modes of integration (...) requires tasks and situations that are not simple adaptation of paper and pencil tasks, often tasks without equivalent in the paper and pencil environment, thus tasks not so easy to design when you enter in the technological world with your paper and pencil culture.  (\cite{icmi}, p. 468)
\end{quote}

The range of hardware and software considered in this Study is huge, from Dynamic Geometry Environments to Computer Algebra Systems, including Animation Microworlds, Games and Spreadsheets, showing that the use of ICT in the mathematics classroom is not limited to any particular kind of software and offers thus many possibilities for mathematics teaching and learning. \\

Another important point visible in this Study is the need to find an answer to the ''wrong-doing'' of
certain technologies. How to deal with the pitfalls of numerical analysis, namely dealing with rounding errors? How to correctly identify a tangent to a circle in a Dynamic Geometry Environment that has difficulties with the continuity? The Study calls for a reasonable 

\begin{quote}
(...) basic understanding of the inner representation of mathematics (e.g., numbers, equations, stochastics, graphical representations, and geometric figures) within a computer and a global awareness of problems related to the difference between conceptual and computational mathematics. (\cite{icmi}, p. 153)
\end{quote}

\section{First Conclusion}

What we discussed from these international studies allows us to conclude that the CCSS are right in investing decidedly in the use of ICT in the classroom:

\begin{quote}
Standards for Mathematical Practice (...) 5 Use appropriate tools strategically.
Mathematically proficient students consider the available tools when solving a
mathematical problem. These tools might include pencil and paper, concrete
models, a ruler, a protractor, a calculator, a spreadsheet, a computer algebra system,
a statistical package, or dynamic geometry software. Proficient students are
sufficiently familiar with tools appropriate for their grade or course to make sound
decisions about when each of these tools might be helpful, recognizing both the
insight to be gained and their limitations. For example, mathematically proficient
high school students analyze graphs of functions and solutions generated using a
graphing calculator. They detect possible errors by strategically using estimation
and other mathematical knowledge. When making mathematical models, they know
that technology can enable them to visualize the results of varying assumptions,
explore consequences, and compare predictions with data. Mathematically
proficient students at various grade levels are able to identify relevant external
mathematical resources, such as digital content located on a website, and use them
to pose or solve problems. They are able to use technological tools to explore and
deepen their understanding of concepts.\cite{mcc}
\end{quote}

Most of the countries in the world have a clear vision of what needs to be done. For example the official curriculum for Singapore reads:

\begin{quote}
AIMS OF MATHEMATICS EDUCATION IN SCHOOLS: (...) (6) Make effective use of a variety of mathematical tools (including information
and communication technology tools) in the learning and application of
mathematics. (...)\\
(...) The use of manipulatives (concrete
materials), practical work, and use of technological aids should be part of the
learning experiences of the students.\\
SKILLS: (...) Skill proficiencies include the ability to use technology confidently, where
appropriate, for exploration and problem solving.\cite{sing}
\end{quote}

What happens in the real classroom is not so simple.

\section{A difficult task}

In a national examination in Portugal for the 12th grade, a mathematical modeling problem involved the study of the function:

\begin{equation}
 d(x)=149.6(1-0.0167 \cos x)
\end{equation}

\begin{figure}
\centering
\includegraphics[height=5cm]{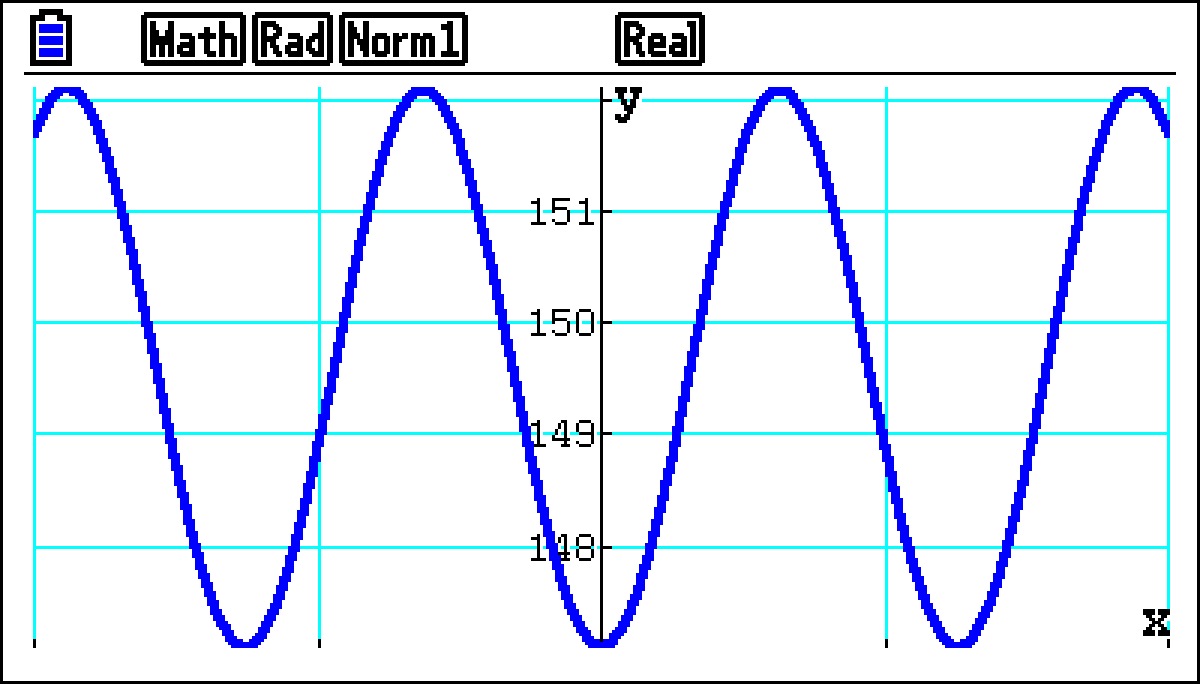} 
\includegraphics[height=5cm]{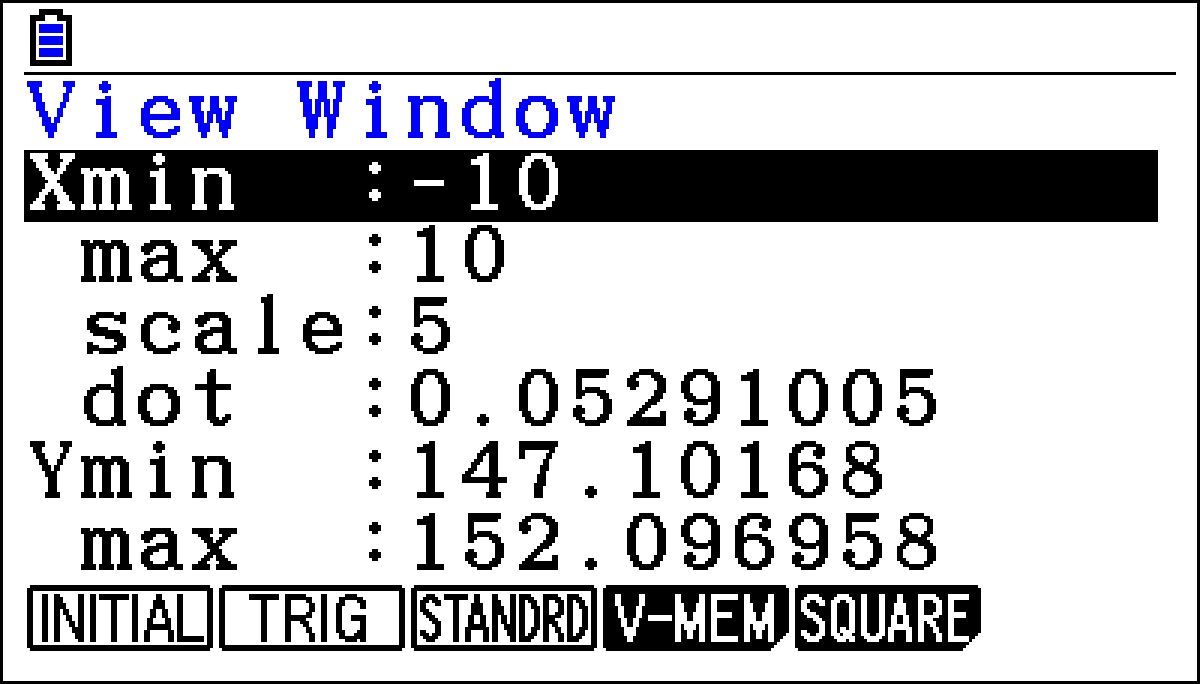}
\caption{Using the AUTO feature to draw a simple graph of $d(x)$ in a GC.}
\end{figure}

Graphing Calculators (GC) are allowed in national examinations in Portugal and so the students can use them to study this function. The biggest challenge here for the student, and it has been proven to be a big obstacle, is to find a viewing window to obtain the graph for this function. Of course you can get the help of the AUTO function of the calculator but then you are not sure you get all the details you need in graph that shows up.\\

The second more difficult group of questions in the national 12th grade examinations in Portugal, were the ones requiring the use of graphing calculators, some of them also involving modeling problems. In a previous study \cite{bal} we concluded that these problems all involve the need to choose a viewing window. There is no algorithm that can guarantee you get the best viewing window. You can produce a table of values to help you but you will need always to experiment or know some properties of the function in order to be sure you get a ''complete'' graph. In more difficult situations you may need to use more than one graph to capture the details of the graph of the function you want to study.\\

Another similar difficulty is discussed by Luc Trouche in his paper \cite{trou} in the journal \emph{Educational Studies in Mathematics}. If a student tries to use a GC to study the limit of a function when the independent variable goes to $+ \infty$ he will try to graph the function ''as far as possible''. But if he is faced with a function like

\begin{equation}
 f(x)=\ln x + 10 \sin x
\end{equation}

he will think it will not have a limit, when the limit is really $+ \infty$. The graph will give him a dangerous message:

\begin{figure}
\centering
\includegraphics[height=5cm]{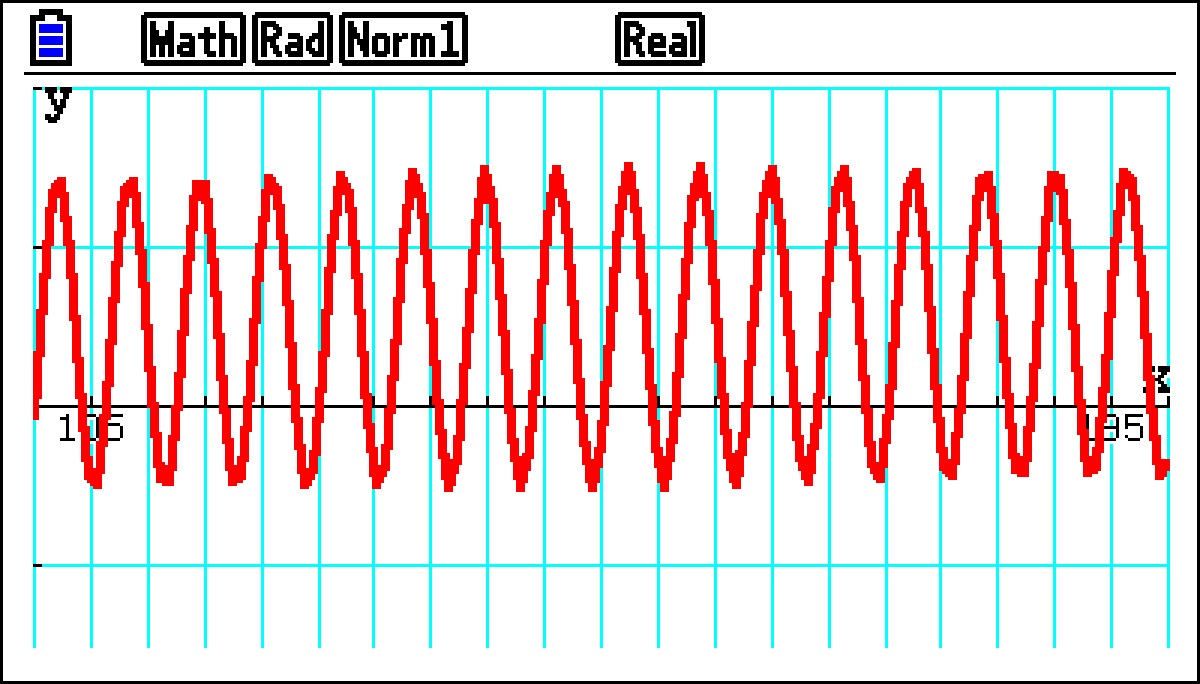} 
\caption{Graphing function $f(x)$ in the viewing window $[100,200] \times [-15,20]$.}
\end{figure}

Luc Trouche concludes that the complexity of the transformation of this new artifact into a useful instrument for the work of the student is related with the sophistication of the artifact, namely when it comes to a GC with CAS-Computer Algebra System. This is a big educational challenge and should be kept in mind when some hardware or software is selected to be used in the classroom.\\

In the same paper Luc Trouche warns against the lack of investment in the use of ICT in the classroom, observing that the use in the classroom is too limited, in France and other countries, and consequently the ''learning of the use of instruments is made most of the time alone or between friends'' (\cite{trou}, p. 190) with all its dangers. It is clear from the international studies that the use of ICT can produce an improvement in student performance, but to arrive at that point of improvement, a lot of research, experimentation and planning must be made. \\

We conclude with a recommendation made by Seymour Papert, the inventor of the educational programming language LOGO, in the plenary talk he gave at the Study Conference included in ICMI Study 17. Seymour Papert ended his talk asking us to spend reasonable part of our time and energy thinking about possible futures, freeing our minds from the current constraints.

\end{document}